\documentclass[a4paper,11pt]{amsart}%
\usepackage{amsfonts}
\usepackage{amssymb}
\usepackage[utf8]{inputenc}
\usepackage{amsmath}
\usepackage{xcolor}
\usepackage{graphicx}
\usepackage{mathptmx}
\providecommand{\U}[1]{\protect\rule{.1in}{.1in}}
\usepackage[colorlinks=true, linkcolor=red, citecolor=blue]{hyperref}
\def\r{\mathbb{R}}

\def\<{\langle}
\def\>{\rangle}

\def\2{L^2}

\theoremstyle{remark}
\newtheorem{remark}{Remark}
\newtheorem*{acknowledgement}{Acknowledgement}

\begin{document}
\title[IBVP of the KdV]{Initial boundary value problem for Korteweg-de Vries equation: a review and open problems}
\author[Capistrano--Filho]{Roberto A. Capistrano--Filho}
\address{Departmento de Matem\'atica,  Universidade Federal de Pernambuco 50740-545, Recife (PE), Brazil.}
\email{capistranofilho@dmat.ufpe.br}
\author[Sun]{Shu-Ming Sun}
\address{Department of Mathematics, Virginia Tech, Blacksburg, VA 24061, United State}
\email{sun@math.vt.edu}
\author[Zhang]{Bing-Yu Zhang}
\address{Department of Mathematical Sciences, University of Cincinnati, Ohio 45221-0025, United State}
\email{zhangb@ucmail.uc.edu}

\begin{abstract}
In the last 40 years the study of initial boundary value problem for the Korteweg-de Vries equation has had the attention of researchers from various research fields. In this note we present a review of the main results about this topic and also introduce interesting open problems which still requires attention from the mathematical point of view.
\end{abstract}
\maketitle

\section{Introduction}

In 1834 John Scott Russell, a Scottish naval engineer, was observing the Union Canal in Scotland when he unexpectedly witnessed a very special physical phenomenon which he called a wave of translation \cite{Russel1844}. He saw a particular wave traveling through this channel without losing its shape or velocity, and was so captivated by this event that he focused his attention on these waves for several years, not only built water wave tanks at his home conducting  practical and theoretical research into these types of waves, but also challenged the mathematical community to prove theoretically the existence of his solitary waves and to give an a priori demonstration a posteriori.


A number of researchers took up Russell’s challenge. Boussinesq  was  the  first  to  explain  the  existence  of  Scott  Russell’s  solitary wave  mathematically.  He  employed  a  variety  of  asymptotically  equivalent equations to describe water waves in the small-amplitude, long-wave regime. In fact, several  works  presented  to  the  Paris  Academy  of  Sciences  in  1871  and 1872, Boussinesq addressed the problem of the persistence of solitary waves of permanent form on a fluid interface \cite{Boussinesq,Boussinesq1,Boussinesq2,Boussinesq3}. 
It is important to mention that in 1876, the English physicist Lord Rayleigh obtained a different result \cite{Rayleigh}.


After Boussinesq theory,  the Dutch mathematicians D. J. Korteweg and his student G. de Vries derived 
a nonlinear partial differential equation in 1895 that  possesses  a solution describing the phenomenon discovered by Russell, 

\begin{equation}\label{kdv}\frac{\partial\eta}{\partial{t}}=\frac{3}{2}\sqrt{\frac{g}{l}}\frac{\partial}{\partial{x}}\left(\frac{1}{2}\eta^2+\frac{3}{2}\alpha\eta+\frac{1}{3}\beta\frac{\partial^2\eta}{\partial{x^2}}\right),
\end{equation}
in which $\eta$ is the surface elevation above the equilibrium level, $l$ is an arbitrary constant related to the motion of the liquid, $g$ is the gravitational constant, and $\beta=\frac{l^3}{3}-\frac{Tl}{\rho g}$ with surface capillary tension $T$ and density $\rho$. 
The equation (\ref{kdv})  is called   the Korteweg-de Vries equation in the literature, often abbreviated as the KdV equation,   although  it had appeared explicitly in Boussinesq’s massive 1877 Memoir \cite{Boussinesq3}, as equation (283bis) in  a  footnote  on  page  360\footnote{The interested readers are referred to \cite{jager2006, pego1998} for history and origins  of the Korteweg-de Vries equation.}.

Eliminating the physical constants by using the following change of variables
$$t\to\frac{1}{2}\sqrt{\frac{g}{l\beta}}t, \quad x\to-\frac{x}{\beta}, \quad u\to-\left(\frac{1}{2}\eta+\frac{1}{3}\alpha\right)$$
one obtains the standard Korteweg-de Vries equation  $$u_t + 6uu_x + u_{xxx}= 0,$$
which is 
now commonly accepted as a mathematical model for the unidirectional propagation of small-amplitude long waves in nonlinear dispersive systems.


This note is concerned with the main results already obtained for the initial-boundary value problem (IBVP) of the KdV equation posed on a finite interval $(0,L)$.  The first paper which treated this problem was given by Bubnov in 1979  \cite{Bubnov79} when he considered the IBVP of the KdV equation on the finite interval $(0,1)$ with general boundary conditions. After that, many authors worked on improving the existing results and presenting new results in the last 30 years. 

Our intention here is to present the main results on this field. Also, we give some further comments and, at the end, discuss open problems related to the IBVP of the KdV equation in a bounded domain.

\section{A review of IBVP for KdV}
Consider the IBVP of the KdV equation posed on a finite interval $(0,L)$
\begin{equation}\label{1.1}
 u_t+u_x+u_{xxx}+uu_x=0,\qquad u(x,0)=\phi(x), \qquad  0<x<L,  \
 t>0
    \end{equation}
with  general   non-homogeneous boundary conditions posed on the two ends of the
interval  $(0,L)$,
\begin{equation} \label{1.2}
 B_1u=h_1(t), \qquad B_2 u= h_2 (t),  \qquad B_3 u= h_3 (t) \qquad
 t>0,
 \end{equation} where
\[ B_i u =\sum _{j=0}^2 \left(a_{ij} \partial ^j_x u(0,t) + b_{ij}
\partial ^j_x u(L,t)\right), \qquad i=1,2,3,\]
and $a_{ij}, \ b_{ij}$, $ j=0, 1,2, \ i=1,2,3,$ are real constants. The following natural question arises:

\smallskip

\emph{Under what assumptions on the coefficients $a_{ij}, \ b_{ij}
$ in (\ref{1.2}), is the IBVP (\ref{1.1})-(\ref{1.2}) well-posed
in the classical Sobolev space $H^s (0,L)$?}

\smallskip

As mentioned before, Bubnov \cite{Bubnov79} studied  the following IBVP of
the KdV equation on the finite interval $(0,1)$:
\begin{equation}\label{1.3}
\begin{cases}
 u_t +uu_x+u_{xxx}=f, \quad u(x,0)=0, \quad x\in (0,1), \ t\in
 (0,T), \\ \alpha _1 u_{xx}(0,t)+\alpha _2 u_x (0,t)+\alpha _3
 u(0,t)=0, \\ \beta_1 u_{xx} (1,t)+\beta _2 u_x (1,t)+ \beta _3
 u(1,t) =0, \\ \chi _1 u_x (1,t)+ \chi _2 u(1,t)=0
 \end{cases}
 \end{equation}
and obtained the following result.

\smallskip

 \noindent
 {\bf Theorem $\mathcal{A}$ }\cite{Bubnov79}:
\emph{Assume that
\begin{equation}\label{1.4}
\begin{cases}
if \ \alpha _1 \beta_1 \chi _1 \ne 0, \ then \ F_1>0, \ F_2 >0, \\
if  \ \beta _1\ne 0, \ \chi _1 \ne 0, \ \alpha _1 =0, \ then \
\alpha _2=0, \ F_2 >0, \ \alpha _3 \ne 0, \\ if \ \beta _1 =0, \
\chi _1 \ne 0, \ \alpha _1 \ne 0, \ then \ F_1 >0, \ F_3 \ne 0, \\
if \ \alpha _1=\beta _1 =0, \ \chi _1 \ne 0, \ then \ F_3\ne 0, \
\alpha _2 =0, \ \alpha _3 \ne 0, \\ if \ \beta _1 =0, \ \alpha _1
\ne 0, \ \chi _1 =0, \ then \ F_1 >0, \ F_3 \ne 0, \\ if \ \alpha
_1=\beta _1 =\chi _1 =0, \ then \ \alpha _2 =0, \ \alpha _3 \ne 0, \
F_3 \ne 0,
\end{cases}
\end{equation}
where
\[ F_1=\frac{\alpha _3}{ \alpha _1} -\frac{\alpha _2^2}{2\alpha
_1^2}, \ F_2 =\frac{\beta_2 \chi _2}{\beta _1 \chi _1}
-\frac{\beta _3}{\beta _1} -\frac{\chi _2^2}{2\chi _1^2}, \ F_3
=\beta _2 \chi_2-\beta _1\chi _1 . \]  For any given
\[ f\in H^1_{loc}(0, \infty ; L^2 (0,1)) \ with \ \ f(x,0)=0, \]
 there exists a $T>0$ such that (\ref{1.3}) admits a unique
solution
\[ u\in L^2 (0, T; H^3 (0,1)) \ with  \ u_t \in L^{\infty} (0,T;
L^2 (0,1))\cap L^2 (0,T; H^1 (0,1)) .\] }

The main tool used by Bubnov to prove his theorem is the
following Kato type smoothing property for solution $u$ of the
linear system associated to the  IBVP (\ref{1.3}),
\begin{equation}\label{1.5}
\begin{cases}
 u_t  +u_{xxx}=f, \quad u(x,0)=0, \quad x\in (0,1), \ t\in
 (0,T), \\ \alpha _1 u_{xx}(0,t)+\alpha _2 u_x (0,t)+\alpha _3
 u(0,t)=0, \\ \beta_1 u_{xx} (1,t)+\beta _2 u_x (1,t)+ \beta _3
 u(1,t) =0, \\ \chi _1 u_x (1,t)+ \chi _2 u(1,t)=0.
 \end{cases}
 \end{equation}
 Under the
 assumptions (\ref{1.4}):
\begin{equation*}f\in L^2(0,T; L^2 (0,1))\implies  u\in
L^2 (0,T; H^1 (0,1))\cap L^{\infty} (0,T; L^2 (0,1))
\end{equation*}
and
\[ \|u\|_{L^2 (0,T; H^1 (0,1))}+ \|u\|_{L^{\infty} (0,T; L^2
(0,1))} \leq C\|f\|_{L^2 (0,T; L^2 (0,T))} \] where $C>0$ is a
constant independent of $f$.

  In the past thirty years since the work of
Bubnov, various boundary-value problems of the KdV equation have
been studied. In particular, the following three classes of IBVPs of the
KdV equation on the finite interval $(0,L)$,
\begin{equation}\label{1.7}
\begin{cases}
u_t +u_x +uu_x +u_{xxx}=0, \ u(x,0)=\phi (x), \quad x\in (0,L), \
t>0, \\ u(0,t)= h_1(t), \quad u(L,t) = h_2 (t), \quad u_x (L,t)
=h_3 (t), \end{cases}
\end{equation}
\begin{equation}\label{1.8}
\begin{cases}
u_t +u_x +uu_x +u_{xxx}=0, \ u(x,0)=\phi (x), \quad x\in (0,L), \
t>0, \\ u(0,t)= h_1(t), \quad u(L,t) = h_2 (t), \quad u_{xx} (L,t)
=h_3 (t), \end{cases}
\end{equation}
and
\begin{equation}\label{1.8a}
\begin{cases}
u_t +u_x +uu_x +u_{xxx}=0, \ u(x,0)=\phi (x), \quad x\in (0,L), \
t>0, \\ u_{xx}(0,t)= h_1(t), \quad u_x(L,t) = h_2 (t), \quad u_{xx} (L,t)
=h_3 (t), \end{cases}
\end{equation}
as well as the IBVPs of the KdV equation posed in a quarter plane have been intensively studied in the past twenty years (cf. \cite{BSZ03FiniteDomain,bsz-finite,ColGhi97,Fam83,Fam89, faminskii2004, faminskii2007, Holmer06,KrZh,KrIvZh,RiUsZh}  and the references therein)
following the rapid advances of the study of the pure initial value
problem of the KdV equation posed on the whole line $\r$ or on the
periodic domain $\mathbb{T}$ (cf. \cite{BS76,BS78,Bourgain93a,Bourgain93b,Bourgain97,ColKeel03,Fam83,Fam89,Fam99,KPV89,KPV91,KPV91-1,KPV93,KPV93b,KPV96}  and the references therein).

The nonhomogeneous IBVP (\ref{1.7}) was first studied by Faminskii in \cite{Fam83,Fam89}
and was shown to be well-posed in the  spaces $L^2 (0,L)$ and $H^3 (0,L)$.

\smallskip

\noindent {\bf Theorem $\mathcal{B}$} \cite{Fam83,Fam89} \emph{Let $T>0$ be given. For
any $\phi \in L^2 (0,L)$	 and $\vec{h}= (h_1, h_2, h_3)$ belonging of
\[ W^{\frac13, 1}(0,T)\cap
L^{6+\epsilon} (0,T)\cap H^{\frac16} (0,T)\times W^{\frac56
+\epsilon, 1} (0,T)\cap H^{\frac13} (0,T)\times  L^2 (0,T),\]
the IBVP (\ref{1.7}) admits a unique solution $$u\in C([0,T]; L^2
(0,L))\cap L^2 (0,T; H^1 (0,L)).$$ Moreover, the solution map is
continuous in the corresponding spaces.}
  \emph{In addition, if $\phi \in
H^3 (0,L)$, $ h_1' \in W^{\frac13, 1}(0,T)\cap L^{6+\epsilon}
(0,T)\cap H^{\frac16} (0,T)$, $h_2'\in W^{\frac56 +\epsilon, 1} (0,T)\cap H^{\frac13 } (0,T)$  and $ h_3' \in L^2 (0,T)$
 with
$$ \phi (0)=h_1 (0), \phi (L)=h_2 (0), \ \phi' (L) = h_3 (0),$$
then the solution $u\in C^1([0,T]; H^3 (0,L))\cap L^2 (0,T; H^4(0,L))$.}

\smallskip

Bona \textit{et al.} in \cite{BSZ03FiniteDomain} showed that the IBVP (\ref{1.7}) is
locally well-posed in the space $H^s (0,L)$ for any $s\geq 0$:

\smallskip
\noindent {\bf Theorem $\mathcal{C}$} \cite{BSZ03FiniteDomain}:
 \emph{Let $s\geq
0$ , $r>0$ and $T>0$ be given.}  \emph{There
exists a $T^*\in (0, T]$ such that for  any $s-$compatible
$\phi \in H^s (0,L)$ and
\[ \vec{h}= (h_1, h_2, h_3) \in
H^{\frac{s+1}{3}}(0,T) \times H^{\frac{s+1}{3}}(0,T)\times
H^{\frac{s}{3}} (0,T) \] satisfying
\[ \| \phi \| _{H^s (0,L)} + \| \vec{h}\|_{H^{\frac{s+1}{3}}(0,T) \times H^{\frac{s+1}{3}}(0,T)\times
H^{\frac{s}{3}} (0,T)} \leq r,\] the IBVP (\ref{1.7}) admits a
unique solution
\[ u\in C([0,T^*]; H^s (0,L))\cap L^2 (0,T^*; H^{s+1} (0,L)).\]
Moreover, the corresponding solution map is  analytic
in the corresponding spaces.}

\smallskip
Holmer \cite{Holmer06} proved that IBVP (\ref{1.7}) is locally
well-posed in the space $H^s (0,L)$ for any $-\frac34 <s<
\frac12$,  and Bona \textit{et al.} in \cite{bsz-finite} showed that the IBVP
(\ref{1.7}) is locally well-posed $H^s (0,L)$ for any $s>-1$.


As for the IBVP (\ref{1.8}), its study began with the work of Colin and Ghidalia in late 1990's \cite{ColGhi97,ColGhi97a,ColGhi01}.
They obtained  in  \cite{ColGhi01}  the following results.
 \begin{itemize} \item[(i)]  \emph{Given
$h_j\in C^1([0, \infty)), \ j=1,2,3$ and $\phi \in H^1 (0,L)$
satisfying $h_1(0)=\phi (0)$, there exists a $T>0$ such that the
IBVP (\ref{1.8}) admits a solution (in the sense of distribution)}
\[ u\in L^{\infty}(0,T; H^1(0,L))\cap C([0,T]; L^2 (0,L)) .\]

\item[(ii)] \emph{The solution $u$ of the IBVP (\ref{1.8}) exists
globally in $H^1(0,L)$ if the size of its initial value $\phi \in
H^1 (0,L)$ and its boundary values $h_j\in C^1([0, \infty )), \
j=1,2,3$ are all small.}
\end{itemize}
In addition, they showed that the associate linear IBVP
\begin{equation}\label{1.9}
        \begin{cases}
        u_t+u_x+u_{xxx}=0,\qquad u(x,0)=\phi(x)  & x\in (0,L),  \ t\in \r^+ \\
        u(0,t)=0,\ u_x(L,x)=0,\ u_{xx}(L,t)=0
        \end{cases}
    \end{equation}
 possesses a strong smoothing property:

\emph{ For any $\phi \in L^2 (0,L)$, the linear IBVP (\ref{1.9})
admits a
 unique solution $$u\in C(\r^+; L^2 (0,L))\cap L^2 _{loc} (\r^+; H^1
 (0,L)).$$}
Aided by this smoothing property, Colin and Ghidaglia  showed
that the homogeneous IBVP (\ref{1.8}) is locally well-posed in the
space $L^2 (0,L)$.

\smallskip

 \noindent
 {\bf Theorem $\mathcal{D}$} \cite{ColGhi01}  \emph{Assuming $h_1=h_2=h_3\equiv 0$, then for any given $\phi \in L^2
(0,L)$, there exists a $T>0$ such that the IBVP (\ref{1.8}) admits
a unique weak solution $u\in C([0,T]; L^2 (0,L))\cap L^2 (0,T; H^1
(0,L))$.}

\smallskip

Returning the attention to the IBVP \eqref{1.8}, Rivas \textit{et al.} in \cite{RiUsZh}, showed that the solutions exist globally as long as their initial values and the associated boundary data are small, they proved the following result:

\smallskip

 \noindent
 {\bf Theorem $\mathcal{E}$} \cite{RiUsZh} \emph{Let $s\geq 0$ with
$s\neq\frac{2j-1}{2}, \text{} \text{}j=1,2,3...$
There exist positive constants $\delta$ and $T$ such that for any $s-$compatible
 $\phi \in H^s (0,L)$ and $\vec{h}= (h_1, h_2, h_3)$ on the class
$$B^s_{(t,t+T)}:=H^{\frac{s+1}{3}}(t,t+T) \times H^{\frac{s}{3}}(t,t+T)\times
H^{\frac{s-1}{3}} (t,t+T) $$
 with $ \|\phi \|_{H^s (0,L)} + \|\vec{h}\|_{B^s_{(t,t+T)}} \leq \delta,$
 and
 $\sup_{t\geq0}\|\vec{h}\|_{B^s_{(t,t+T)}}<\infty,$
 the IBVP (\ref{1.9}) admits a unique solution
 \[ u\in Y^s_{(t,t+T)}:=C([t,t+T]; H^s (0,L))\cap L^2 (t,t+T; H^{s+1}(0,L))\] such that for any $t\geq0,$ $\sup_{t\geq0}\|\vec{v}\|_{Y^s_{(t,t+T)}}<\infty.$}


\smallskip

More recently, Kramer \textit{et al.} in \cite{KrIvZh} showed that the IBVP \eqref{1.8} is locally well-posedness in the classical Sobolev space $H^s(0,L)$, for $s>-\frac{3}{4}$, which provides a positive answer to one of the open questions of Colin and Ghidalia \cite{ColGhi01}.

 Kramer and Zhang in \cite{KrZh}, studied the following non-homogeneous
boundary value problem,
\begin{equation}\label{1.3-g}
\begin{cases}
 u_t +uu_x+u_{xxx}=0, \quad u(x,0)=\phi (x), \quad x\in (0,1), \ t\in
 (0,T), \\ \alpha _1 u_{xx}(0,t)+\alpha _2 u_x (0,t)+\alpha _3
 u(0,t)=h_1(t), \\ \beta_1 u_{xx} (1,t)+\beta _2 u_x (1,t)+ \beta _3
 u(1,t) =h_2(t), \\ \chi _1 u_x (1,t)+ \chi _2 u(1,t)=h_3 (t).
 \end{cases}
 \end{equation}
 They showed that the IBVP (\ref{1.3-g}) is locally well-posed in
 the space $H^s (0,1)$ for any $s\geq 0$ under the assumption (\ref{1.4}).

\smallskip

 \noindent
 {\bf Theorem $\mathcal{F}$} \cite{KrZh} \emph{Let $s\geq 0$ and $T>0$ be given and assume (\ref{1.4}) holds. For any $r>0$,
 there exists a $T^*\in (0,T]$ such that for any $s-$compatible
$\phi \in H^s (0,1)$, $h_j\in H^{\frac{s+1}{3}}(0,T),
 j=1,2,3$ with
 \[ \|\phi \|_{H^s (0,1)} + \|h_1\|_{H^{\frac{s+1}{3}}(0,T)}
 +\|h_2\|_{H^{\frac{s+1}{3}}(0,T)}+\|h_3\|_{H^{\frac{s+1}{3}}(0,T)}
 \leq r,\]  the IBVP (\ref{1.3-g}) admits a unique solution
 \[ u\in C([0,T^*]; H^s (0,1))\cap L^2 (0,T^*; H^{s+1}(0,1)) .\]
 Moreover, the solution $u$ depends continuously on its initial data
 $\phi $ and the boundary values $h_j, j=1,2,3$ in the respective
 spaces.}

Recently, Capistrano--Filho \textit{et al.} \cite{CCFZh} studied the IBVP \eqref{1.8a}. The authors proved the local well-posedness for this system. More precisely:

\smallskip

 \noindent {\bf Theorem $\mathcal{G}$} \cite{CCFZh} \emph{ Let $T>0$ and $s\geq0$. There exists a $T^*\in(0,T]$ such that for any $(\phi , \vec{h}) \in X_{
T}$, where
\[ X_{T}:= H^s  (0,L)\times H^{\frac{s-1}{3}}(0,T)\times
H^{\frac{s}{3}}(0,T)\times H^{\frac{s-1}{3}}(0,T)\]
 the IBVP \eqref{1.8a} admits a unique solution
$$
 u\in C([0,T];H^s (0,L))\cap L^2(0,T;H^{s+1}(0,L))
$$
In addition, the solution $u$  possesses the hidden regularities
$$\partial_x^lu\in L^{\infty}(0,L;H^{\frac{s+1-l}{3}}(0,T^*)) \quad \text{ for }\quad l=0,1,2.$$
and, moreover, the corresponding solution map is Lipschitz continuous.}

Finally, in a recently work, Capistrano--Filho \textit{et al.}  in \cite{CaSunZha2018} studied the well-posedness of IBVP (\ref{1.1})-(\ref{1.2}). The authors proposed the following hypotheses on those coefficients $a_{ij},
\ b_{ij}$, $ j,i=0, 1,2,3$:
\begin{itemize}
\item[(A1)] $ a_{12}=a_{11}=0, \ a_{10}\ne0, \
b_{12}=b_{11}=b_{10}=0$;

\item[(A2)] $a_{12}\ne0, \ b_{12}=0$;

\item[(B1)] $b_{22}=b_{21}=0,\ b_{20}\ne0, \ a_{22}=a_{21}=a_{20}
=0$;

\item[(B2)] $b_{22}\ne 0, \ a_{22}=0 $;

\item[(C)] $b_{32}=0, \ b_{31}\ne 0, \ a_{32}=a_{31}=0.$
\end{itemize}
For $s\geq 0$, consider  the set $$H^s_0(0,L):=\{\phi(x)\in H^s(0,L): \phi^{(k)}(0)=\phi^{(k)}(L)=0\}$$
with $k=0,1,2,  \cdots ,  [s]$ and
$$H^s_0(0,T]:=\{h(t)\in H^s(0,T):h^{(j)}(0)=0\},$$
for $j=0,1,...,, [s] $. In addition, letting
\begin{equation*}
\begin{cases}
 {\mathcal{ H}}_1^s (0,T) := H_0^{\frac{s+1}{3}}(0,T]\times
H_0^{\frac{s+1}{3}}(0,T]\times H_0^{\frac{s}{3}}(0,T], \\ {\mathcal{ H}}^s_2 (0,T):= H_0^{\frac{s+1}{3}}(0,T]\times
H_0^{\frac{s-1}{3}}(0,T]\times H_0^{\frac{s}{3}}(0,T],\\
{\mathcal{ H}}^s_3 (0,T):= H_0^{\frac{s-1}{3}}(0,T]\times
H_0^{\frac{s+1}{3}}(0,T]\times H_0^{\frac{s}{3}}(0,T], \\ {\mathcal{ H}}^s_4 (0,T):= H_0^{\frac{s-1}{3}}(0,T]\times
H_0^{\frac{s-1}{3}}(0,T]\times H_0^{\frac{s}{3}}(0,T]
\end{cases}
\end{equation*}
and
\begin{equation*}
\begin{cases}
 {\mathcal{ W}}_1^s (0,T) := H^{\frac{s+1}{3}}(0,T)\times
H^{\frac{s+1}{3}}(0,T)\times H^{\frac{s}{3}}(0,T), \\ {\mathcal{ W}}^s_2 (0,T):= H^{\frac{s+1}{3}}(0,T)\times
H^{\frac{s-1}{3}}(0,T)\times H^{\frac{s}{3}}(0,T),\\
{\mathcal{ W}}^s_3 (0,T):= H^{\frac{s-1}{3}}(0,T)\times
H^{\frac{s+1}{3}}(0,T)\times H^{\frac{s}{3}}(0,T), \\ {\mathcal{ W} }^s_4 (0,T):= H^{\frac{s-1}{3}}(0,T)\times
H^{\frac{s-1}{3}}(0,T)\times H^{\frac{s}{3}}(0,T),
\end{cases}
\end{equation*}
they proved the following well-posedness results for the  IBVP (\ref{1.1})-(\ref{1.2}):

\smallskip

 \noindent
 {\bf Theorem $\mathcal{H}$} \cite{CaSunZha2018} \emph{ Let $s\geq 0$ with
$s\neq\frac{2j-1}{2}, \text{} \text{}j=1,2,3..., $ and
$T>0$ be given.  If one  of the assumptions below is satisfied,
\begin{itemize}
\item[(i)] (A1), (B1) and (C) hold,
\item[(ii)] (A1), (B2) and (C) hold,
\item[(iii)] (A2), (B1) and (C) hold,
\item[(iv)] (A2), (B2) and (C) hold,
\end{itemize}
then,  for any $r>0$,  there exists a $T^*\in (0, T]$ such that for any $$(\phi , \vec{h})\in H^s_0
(0,L)\times{\mathcal H}^s_1(0,T)$$ satisfying $\|(\phi, \vec{h})\|_{L^2 (0,L)\times{\mathcal H}^0_1(0,T)} \leq r$,  the IBVP (\ref{1.1})-(\ref{1.2})
admits a solution
$$u\in C([0,T^*]; H^s (0,L))\cap L^2 (0,T^*;H^{s+1}(0,L))$$
possessing the hidden regularity (the sharp Kato smoothing properties)
$$\partial_x^lu\in L^{\infty}(0,L;H^{\frac{s+1-l}{3}}(0,T^*)) \quad \text{ for }\quad l=0,1,2.$$
Moreover, the corresponding solution map is analytically  continuous.}

\section{Further comments}

Before presenting the main ideas to prove Theorem $\mathcal{G}$, let us introduce the following boundary operators $\mathcal{B}_k,  \ k=1,2,3,4$ as $\mathcal{B}_k = \mathcal{B}_{k,0}+ \mathcal{B}_{k,1}$
with
 \[
      \mathcal{B}_{1,0}v:= ( v(0,t), v (L,t), v_{x}(L,t)), \quad
\mathcal{B}_{2,0}v:=( v(0,t), v_{x}(L,t),
        v_{xx}(L,t) ),\]
\[
      \mathcal{B}_{3,0}v:=( v_{xx}(0,t), v(L,t),   v_{x}(L,t)), \quad   \mathcal{B}_{4,0}v:=( v_{xx}(0,t), v_{x}(L,t),
    v_{xx}(L,t))\]
 and
\begin{align*} 
& \mathcal{B}_{1,1}v:=\left (0,\ 0, \ 0\right ),\\
&\mathcal{B}_{2,1}v:=\left (0,\  b_{30}v(L,t),\ a_{21} v_x (0,t) + b_{20} v(L,t)\right),\\
&\mathcal{B}_{3,1}v:=\left(a_{10}  v(0,t) +a_{11} v_x (0,t),\ 0,\ a_{30}v(0,t) \right),\\
& \mathcal{B}_{4,1}v:=\left( \sum _{j=0}^1 a_{1j} \partial ^j_x v(0,t)+ b_{10} v(L,t),\  a_{30}v(0,t)+b_{30}v(L,t), \  \sum _{j=0}^1 a_{2j} \partial ^j_x v(0,t) +  b_{20} v(L,t)\right).
\end{align*}
Thus, the  assumptions imposed on the boundary conditions in Theorem $\mathcal{G}$ can be reformulated  as follows:
\begin{itemize}
\item[(i)]  $((A1), (B1), (C)) \Leftrightarrow \mathcal{B}_{1}v= \vec{h},$
\item[(ii)]  $((A1),(C),(B2)) \Leftrightarrow \mathcal{B}_{2}v= \vec{h},$
\item[(iii)]  $((A2), (B1), (C)) \Leftrightarrow \mathcal{B}_{3}v= \vec{h},$
\item[(iv)]  $((A2),(C),(B2)) \Leftrightarrow \mathcal{B}_{4}v= \vec{h}.$
\end{itemize}
In \cite{CaSunZha2018}, to prove Theorem $\mathcal{G}$, the authors first studied the linear IBVP
\begin{equation}\label{y-1}
\begin{cases}
u_t +u_{xxx} +\delta_k u=f, \quad x\in (0,L), \quad t >0\\ u(x,0)= \phi (x), \\ \mathcal{B}_{k,0} u= \vec{h},
\end{cases}
\end{equation}
for $k=1,2,3,4$, to establish all the linear estimates needed for  dealing with the nonlinear IBVP (\ref{1.1})-(\ref{1.2}).  Here $\delta _k=0$ for $k=1,2,3$ and $\delta _4=1$.

After that, they considered the nonlinear map $\Gamma $ defined by the following IBVP:
\begin{equation}\label{y-2}
\begin{cases}
u_t +u_{xxx} +\delta_ku= -v_x -vv_x +\delta_kv , \quad x\in (0,L), \quad t >0\\u(x,0)= \phi (x), \\ \mathcal{B}_{k,0} u= \vec{h}-\mathcal{B}_{k,1} v,
\end{cases}
\end{equation}
showing thus that $\Gamma$ is a contraction in an appropriate space whose  fixed point will be the desired solution of the nonlinear  IBVP (\ref{1.1})-(\ref{1.2}) by using the sharp Kato smoothing  property of the solution of the IBVP (\ref{y-1}).

The main point here is to demonstrate the smoothing properties for solutions of the IBVP (\ref{y-1}).   In order to overcome this difficulty, Capistrano--Filho \textit{et al.} in \cite{CaSunZha2018} needed to study the following IBVP
\begin{equation}\label{y-3}
\begin{cases}
u_t +u_{xxx}+\delta _ku=0, \quad x\in (0,L), \quad t >0,\\ u(x,0)= 0, \\ \mathcal{B}_{k,0} u= \vec{h}.
\end{cases}
\end{equation}
The corresponding solution map $\vec{h} \to u$ will be called the \textit{boundary integral operator} denoted by ${\mathcal W}_{bdr} ^{(k)}$.  An explicit representation  formula is given for this boundary integral operator that plays an important role  in showing the solution of the IBVP (\ref{y-3}) possesses the  smoothing properties. The  needed  smoothing properties for solutions of the IBVP (\ref{y-1}) will then follow from the  smoothing properties for solutions of the IBVP (\ref{y-3}) and  the well-known sharp Kato smoothing properties for solutions of the Cauchy problem
\[ u_t +u_{xxx} +\delta _ku=0, \quad u(x,0)=\psi (x), \quad x, \ t\in \mathbb{R}.\]

Finally, the following comments are now given in order:

\begin{remark} The  temporal regularity conditions imposed on the boundary values  $\vec{h}$ on Theorem $\mathcal{G}$ are optimal (cf. \cite{BSZ02,BSZ04,BSZ06}).
\end{remark}
\begin{remark} As a comparison, note that the assumptions of Theorem $\mathcal{A}$ are equivalent to  one of the following boundary conditions imposed on the equation in (\ref{1.3}):

a)
$$u(0,t)=0, \quad u(1,t)=0, \quad u_x (1,t)=0;$$

b)
$$u_{xx}(0,t)+au_x(0,t)+bu(0,t)=0, \quad u_x(1,t)=0,
 \quad u(1,t)=0$$
 with
 \begin{equation}\label{z-1} a>b^2/2;\end{equation}

c)$$u(0,t)=0, \quad u_{xx}(1,t)+au_x(1,t)+bu(1,t)=0,
 \quad u_x(1,t)+cu(1,t)=0,$$
 with
 \begin{equation}\label{z-2} ac>b-c^2/2;\end{equation}

d) $$u_{xx}(0,t)+a_1u_x(0,t)+a_2 u(0,t)=0, $$ $$ u_{xx}(1,t)+b_1u_x(1,t)+b_2 u(1,t)=0,$$ and $$u_x(1,t)+cu(1,t)=0,$$
 with
 \begin{equation} \label{z-3} a_2 > a_1^2/2, \quad b_1c > b_2 -c^2/2 .\end{equation}
 It follows from Theorem $\mathcal{G}$ that conditions (\ref{z-1}), (\ref{z-2}) and (\ref{z-3}) for Theorem $\mathcal{A}$ can be removed.
 \end{remark}

\section{Open problems}
While the results reported in this paper gave a significant improvement in the theory of initial boundary value problems of the KdV equation on a finite interval,  there are still many questions to be addressed for the following IBVP:
 \begin{equation}\label{4.1}
\begin{cases}
 u_t+u_x+u_{xxx}+uu_x=0, \qquad  0<x<L, \
 t>0,\\
 u(x,0)=\phi(x),\\
{\mathcal B}_ku=\vec{h}.
 \end{cases}
 \end{equation}
 Here we list a few of them which are most interesting to us.


\smallskip
\noindent$\bullet$ {\em Is the IBVP (\ref{4.1}) globally well-posed in the space $H^s (0,L)$ for some $s\geq 0$ or equivalently,  does  any solution of the IBVP (\ref{4.1})  blow up in the some space $H^s (0,L)$ in finite time?}

 \smallskip
 It is not clear  if the IBVP (\ref{4.1}) is globally well-posed  or not even in the  case of $\vec{h}\equiv 0 $.   It follows  Theorem $\mathcal{G}$ (see \cite{CaSunZha2018}) that a solution $u$ of the IBVP (\ref{4.1}) blows up in  the space $H^s (0,L)$   for some $s\geq 0$ at  a finite time $T>0$ if and only   if
 \[ \lim_{t\to T^-}\| u(\cdot, t) \|_{L^2 (0,L)} =+\infty .\]
 Consequently, it suffices to establish a global a priori  $L^2 (0,L)$ estimate
 \begin{equation}\label{priori} \sup _{0\leq t< \infty} \|u(\cdot, t)\|_{L^2 (0,L)} < +\infty ,\end{equation}
 for solutions of the IBVP (\ref{4.1}) in order to obtain the global well-posedness of the IBVP (\ref{4.1}) in the space $H^s (0,L)$ for any $s\geq 0$. However, estimate (\ref{priori}) is known to be held only in one case
 \[ \begin{cases}
 u_t+u_x +uu_x + u_{xxx}=f, \qquad  0<x<L, \
 t>0\\
 u(x,0)=\phi(x)\\
u(0,t)=h_1 (t), \ u(L,t)= h_2 (t), \ u_x (L,t) =h_3 (t). \end{cases}
\]
\smallskip
\noindent$\bullet$ {\em Is the IBVP well-posed in the space $H^s (0,L)$  for some $s\leq -1$?}

\smallskip
Theorem $\mathcal{G}$ ensures that the IBVP (\ref{4.1}) is locally well-posed in the space $H^s (0,L)$ for any $s\geq 0$. Theorem $\mathcal{G}$ can also be extended to the case of $-1< s\leq 0$ using the same approach developed in \cite{bsz-finite}.   For the pure initial value problems  (IVP) of the KdV equation posed on the whole line $\r$ or on torus $\mathbb{T}$,
\begin{equation}
\label{p-1}
u_t +uu_x +u_{xxx}=0, \quad u(x,0)= \phi (x), \quad x, \ t\in \r
\end{equation}
and
\begin{equation}
\label{p-2}
u_t +uu_x +u_{xxx}=0, \quad u(x,0)= \phi (x), \quad x \in \mathbb{T}, \ t\in \r ,
\end{equation}
it is well-known that the IVP (\ref{p-1}) is  well-posed in the space $H^s (\r)$ for any $s\geq -\frac34$ and is  (conditionally) ill-posed in the
space $H^s (\r) $ for any $s< -\frac34$ in the sense the corresponding solution map cannot be uniformly continuous.  As for the IVP (\ref{p-2}), it is well-posed in the space $H^s (\mathbb{T}) $ for any $s\geq -1$.   The solution map corresponding to the IVP (\ref{p-2}) is real analytic when $s>-\frac12$, but only continuous (not even  locally uniformly continuous) when $-1\leq s<-\frac12$. Whether  the IVP (\ref{p-1})  is well-posed  in the space $H^s(\r)$ for any $s<-\frac34$ or the IVP (\ref{p-2}) is well-posed in the space $H^s (\mathbb{T})$  for any $s< -1$ is still an open question. On the other hand, by contrast, the IVP of the KdV-Burgers equation
\[ u_t +uu_x +u_{xxx}-u_{xx}=0, \quad u(x,0)=\phi (x), \quad x\in \r, \ t>0 \]
is known to  be well-posed in the  space $H^s(\r) $ for any $s\geq -1$,  but  is known to be ill-posed for any $s<-1$. We conjecture that the IBVP (\ref{4.1}) is ill-posed in the space $H^s (0,L)$ for any $s<-1$.

\vspace{0.2cm}

Finally, still concerning with well-posedness problem, while the approach developed  recently in \cite{CaSunZha2018} studies the nonhomogeneous boundary value problems of  the KdV equation on $(0,L)$ with quite general boundary conditions,  there are still  some boundary value  problems of the KdV equation that the approach do not work, for example
\begin{equation}
\label{p-3}
\begin{cases}
u_t +uu_x +u_{xxx}=0, \quad x\in (0,L)\\
  u(x,0)= \phi (x), \\
  u(0,t)=u(L,t), \ u_x (0,t)=u_x (L,t), \ u_{xx} (0,t)= u_{xx}(L,t)
  \end{cases}
\end{equation}
and
\begin{equation}
\label{p-4}
\begin{cases}
u_t +uu_x +u_{xxx}=0, \quad x\in (0,L),\\
  u(x,0)= \phi (x),  \\
  u(0,t)=0, \ u (L,t)=0, \ u_{x} (0,t)= u_{x}(L,t) .
  \end{cases}
\end{equation}
A common feature for these two boundary value problems is that the $L^2-$norm of their solutions are conserved:
\[ \int ^L_0 u^2 (x,t) dx =\int ^L_0 \phi ^2 (x) dx \qquad \mbox{for any $t\in \r$}.\]
The IBVP (\ref{p-3}) is equivalent  to the IVP (\ref{p-2}) which was  shown by Kato \cite{Kato79,Kato83} to be well-posed in the space $H^s (\mathbb{T})$ when $s>\frac32$  as early as in the late 1970s.    Its well-posedness in the space $H^s (\mathbb{T})$ when $s\leq \frac32$ , however, was established 24  years later in the celebrated work of Bourgain \cite{Bourgain93a,Bourgain93b} in 1993. As for the IBVP (\ref{p-4}),  its associated linear problem
\begin{equation*}
\begin{cases}
u_t  +u_{xxx}=0, \quad x\in (0,L),\\
  u(x,0)= \phi (x),
   u(0,t)=0, \\ u (L,t)=0, \ u_{x} (0,t)= u_{x}(L,t)
  \end{cases}
\end{equation*}
has been shown by Cerpa  (see, for instance,  \cite{cerpatut})  to be well-posed  in the space $H^s (0,L)$ forward and backward in time. However, the following problem is still unknown:

\smallskip
\noindent$\bullet$ {\em Is the nonlinear IBVP (\ref{p-4}) well-posed in the space $H^s (0,L)$ for some $s$} ?

\subsection{Control theory}
Control theory for KdV equation has been extensively studied in the past two decades and the interested reader is referred to \cite{cerpatut} for an overall view of the subject. As it is possible to see in the paper above, several authors have addressed the study of control theory of the IBVP (see, e.g, \cite{CCFZh,CerRivZha,Rosier97}), who worked on the following four problems related to the IBVP \eqref{4.1}:

\begin{equation*}\label{4.4}
      \mathcal{B}_{1,0}v:=\begin{cases}
          \smallskip u(0,t)=h_{1,1}(t),& t\geq 0,\\
          \smallskip u (L,t)=h_{2,1}(t), & t\geq 0, \\
           u_{x}(L,t)=h_{3,1}(t),& t\geq 0,
        \end{cases}
\quad    \mathcal{B}_{2,0}v:=\begin{cases}
     \smallskip u(0,t)=h_{1,2}(t),& t\geq 0, \\
      \smallskip u_{x}(L,t)=h_{2,2}(t),& t\geq 0, \\
        u_{xx}(L,t)=h_{3,2}(t),& t\geq 0,
        \end{cases}
    \end{equation*}

        \smallskip

\begin{equation*}\label{4.6}
      \mathcal{B}_{3,0}u:=\begin{cases}
      \smallskip u_{xx}(0,t)=h_{1,3}(t)&t\geq 0,\\
      \smallskip u(L,t)=h_{2,3}(t), &t\geq 0,\\
      u_{x}(L,t)=h_{3,3}(t),&t\geq 0
        \end{cases}
   \quad
  \mathcal{B}_{4,0}u:=\begin{cases}
  \smallskip u_{xx}(0,t)=h_{1,4}(t), & t\geq 0,\\
   \smallskip u_{x}(L,t)=h_{2,4}(t), & t\geq 0, \\
    u_{xx}(L,t)=h_{3,4}(t),& t\geq 0.
        \end{cases}
    \end{equation*}

The first class of problem \eqref{4.1}--$\mathcal{B}_{1,0}v$ was studied by Rosier \cite{Rosier97} considering only the control input $h_{3,1}$  (i.e. $h_{1,1}=h_{2,1}=0$). It was shown in \cite{Rosier97} that the exact controllability of the linearized system holds in $L^2(0,L)$ if and only if, $L$ does not belong to the following countable set of critical lengths
\begin{equation*}
\mathcal{N}:=\left\{  \frac{2\pi}{\sqrt{3}}\sqrt{k^{2}+kl+l^{2}}
\,:k,\,l\,\in\mathbb{N}^{\ast}\right\}.
\end{equation*}
The analysis developed in \cite{Rosier97} shows that when the linearized system is controllable, the same is true for the nonlinear case. Note that the converse is false, as it was proved in \cite{cerpa,cerpa1,coron}, that is, the (nonlinear) KdV equation is controllable even when $L$ is a critical length but the linearized system is non controllable.

The existence of a discrete set of critical lengths for which the exact controllability of the linearized equation fails was also noticed by Glass and Guerrero in \cite{GG1} when $h_{2,1}$ is taken as control input (i.e. $h_{1,1}=h_{3,1}=0$). Finally, it is worth mentioning the result by Rosier \cite{Rosier2} and Glass and Guerrero \cite{GG} for which $h_{1,1}$ is taken as control input (i.e. $h_{2,1}=h_{3,1}=0$). They proved that system \eqref{4.1} with boundary conditions $\mathcal{B}_{1,0}v$ is then null controllable, but not exactly controllable, because of the strong smoothing effect.

Recently, Cerpa \textit{et al.} in \cite{CerRivZha} proved similar results to those obtained by Rosier \cite{Rosier97} for the system \eqref{4.1} with boundary conditions $\mathcal{B}_{2,0}v$. More precisely, the authors consider the system with one, two or three controls. In addition, using the well-posedness properties proved by Kramer \textit{et al.} in  \cite{KrIvZh}, they also proved that the controls $h_{i,2}$, $i=1,2,3$ belong to sharp spaces and the locally exact controllability of the linear system associated to \eqref{4.1} holds if, and only if, L does not belong to the following countable set of critical lengths
\begin{equation}
\mathcal{F}:=\left\{ L\in\mathbb{R}^+: L^2=-(a^2+ab+b^2) \text{ with } a,b\in\mathbb{C} \text{ satisfying } \frac{e^a}{a^2}=\frac{e^b}{b^2}=\frac{e^{-(a+b))}}{(a+b)^2}\right\}.
 \label{critical1}
\end{equation}
Moreover, they showed that the nonlinear system \eqref{4.1} with boundary conditions $\mathcal{B}_{2,0}v$ is locally exactly controllable \textit{via} the contraction mapping principle.

Recently, Caicedo \textit{et al.}, in \cite{CCFZh}, proved the controllability results for the system \eqref{1.8a}, that is, system \eqref{4.1} with boundary conditions $\mathcal{B}_{4,0}v$. Naturally, they used the same approaches that have worked effectively for system \eqref{4.1} with boundary condition $\mathcal{B}_{1,0}v$ and $\mathcal{B}_{2,0}v$. In particular, if only $h_{2,4}(t)$ is used, they showed that the system \eqref{4.1} with boundary conditions $\mathcal{B}_{4,0}v$ is \textit{locally exactly controllable} as long as
\begin{equation}
L\notin\mathcal{R}:=\mathcal{N}\cup\left\{k\pi:k\in\mathbb{N}^{\ast}\right\}.
\label{critical_new}
\end{equation}
Thus, with respect of the control issue, a natural and interesting open problem arises here:

\smallskip
\noindent$\bullet$ {\em Is the IBVP \eqref{4.1}, with general boundary condition, controllable?}

\begin{acknowledgement}
The authors wish to thank the referee for his/her valuable comments which improved this paper. Roberto Capistrano-Filho was supported by CNPq (Brazilian Technology Ministry), Project PDE, grants 306475/2017-0, 408181/2018-4 and partially supported by CAPES (Brazilian Education Ministry) and Bing-Yu Zhang was partially supported by NSF of China (11571244, 11231007). 
\end{acknowledgement}

\end{document}